\pdfoutput=1

\documentclass[11pt]{article}
\usepackage[spanish,english]{babel}
\usepackage{amsthm}
\usepackage{epsfig, subfigure}
\usepackage{amscd,latexsym,amssymb}
\usepackage{graphicx}
\usepackage{amsmath,amsfonts,latexsym}
\usepackage{epsfig,graphicx}
\usepackage{url}

\begin{document}

\newtheorem{theorem}{Theorem}
\newtheorem{proposition}{Proposition}
\newtheorem{corollary}{Corollary}
\newtheorem{lemma}{Lemma}

\selectlanguage{english}


\title{
 On the Metric Dimension of Infinite Graphs
 \thanks{Research partially supported by projects JA FQM305, MTM2005-08441-C02-01,
 MEC-MTM2006-01267, DURSI2005SGR00692, MTM2005-08990-C02-01 }
}

\author{
J. C\'aceres\thanks{Universidad de Almer\'{i}a, Almer\'{i}a, Spain, jcaceres@ual.es}
\and C. Hernando\thanks{Universitat Polit\`ecnica de Catalunya, Barcelona, Spain, carmen.hernando@upc.edu}
\and M. Mora\thanks{Universitat Polit\`ecnica de Catalunya, Barcelona, Spain, merce.mora@upc.edu}
\and I. M. Pelayo\thanks{Universitat Polit\`ecnica de Catalunya, Barcelona, Spain, ignacio.m.pelayo@upc.edu}
\and M. L. Puertas\thanks{Universidad de Almer\'{i}a, Almer\'{i}a, Spain, mpuertas@ual.es}
}%



\maketitle





\begin{abstract}
A set of vertices $S$ \emph{resolves} a graph $G$ if every vertex
is uniquely determined by its vector of distances to the vertices
in $S$. The \emph{metric dimension} of a graph $G$ is the minimum
cardinality of a resolving set. In this paper we study the metric
dimension of infinite graphs such that all its vertices have
finite degree. We give necessary conditions for those graphs to
have finite metric dimension and characterize infinite trees with
finite metric dimension. We also establish some results about the
metric dimension of the cartesian product of finite and infinite
graphs, and give the metric dimension of the cartesian product of
several families of graphs.

\vspace{0.5cm}\noindent {\it Key words:} ~infinite graph, locally
finite graph, resolving set, metric dimension, cartesian product

\end{abstract}


\section{Introduction}


Throughout this paper a \emph{graph} $G$ is an ordered pair of disjoint sets $(V,E)$
where $V$ is nonempty and $E$ is a subset of unordered pairs of
$V$. The \emph{vertices} and \emph{edges} of $G$ are the elements
of $V=V(G)$ and $E=E(G)$ respectively.
We say that a graph
$G$ is \emph{finite} (resp. \emph{infinite}) if the set $V(G)$ is
finite (resp. {infinite}). The degree of a vertex $u\in V(G)$ is
the number of edges containing $u$ and is denoted by $\deg_G (u)$,
or simply $\deg (u)$ if the graph $G$ is clear. An  graph
is \emph{locally finite} if every vertex has finite degree.
All the graphs considered in this paper are connected and locally finite.

Infinite graphs have been studied by many authors as D.
K\"onig~\cite{konig}, C. St. J. A. Nash-Williams~\cite{crispin},
C. Thomassen ~\cite{tomas},
 R. Diestel
 and many other authors~\cite{diestel}.



We denote the distance between two vertices $u$ and $v$
by $d_G(u,v)$, or simply $d(u,v)$ if the graph $G$ is clear.
A vertex $x$ in a graph $G$ \emph{resolves} two
vertices $u,v$  if $d(u,x)\neq d(v,x)$. A subset of vertices $S$
is a \emph{a resolving set} of $G$ if for any two vertices, there
exists a vertex in $S$ that resolves them. A resolving set with
minimum cardinality is a \emph{metric basis}. If a graph $G$ has
at least a finite resolving set, the \emph{metric dimension}
$\beta(G)$ of $G$ is the cardinality of a metric basis, otherwise
we say that the metric dimension of $G$ is infinite. If $S=\{
x_1,\dots ,x_n \}$ is a finite set of vertices of $G$, we denote
by $r(u|S)$ the vector of distances from $u$ to the vertices of
$S$, that is, $r(u|S)=(d(u,x_1),\dots ,d(u,x_n))$.  Then, $S$ is a
resolving set if and only if $r(u|S)\not= r(v|S)$ for all vertices
$u\not= v$. If $S$ is a finite metric basis we say that $r(u|S)$
are the \emph{metric coordinates} of vertex $u$ respect to $S$.
A resolving set does not necessarily contain a metric basis. For
example, any two distinct vertices of a
path form a resolving set, but a metric basis of this graph
contains exactly one of its endpoints. This fact
makes more difficult to determine the metric dimension of a graph.

Resolving sets in general graphs were defined by Harary and
Melter~\cite{hararymelter} and Slater~\cite{slater}. Resolving
sets have been widely investigated~\cite{ChaZa,ours,ours2}  and
appear in many diverse areas including coin weighting~\cite{coin},
network discovery and verification~\cite{berliova}, robot
navigation~\cite{landmarks,shanmukha}, connected joins in
graphs~\cite{sebotannier}, and strategies for the Mastermind
game~\cite{mastermind}.


In this paper we study resolving sets in infinite graphs. The
first question that arises is to determine infinite graphs with
finite resolving sets. In Section~\ref{sec.md} we give necessary
conditions for this and in Section~\ref{sec.trees} we characterize
acyclic connected infinite graphs, i.e. \emph{infinite trees},
with finite metric dimension. On the other hand, interesting
examples of graphs can be obtained as cartesian products of
graphs. The metric dimension of cartesian products of finite
graphs is studied in~\cite{ours}. In Section~\ref{sec.cartProd} we
give bounds for the metric dimension of cartesian products of
finite and infinite graphs and determine the metric dimension for
certain families of graphs.


\section{Finite Metric Dimension}\label{sec.md}

\begin{figure}[ht]
\begin{center}
\includegraphics[width=0.8\textwidth]{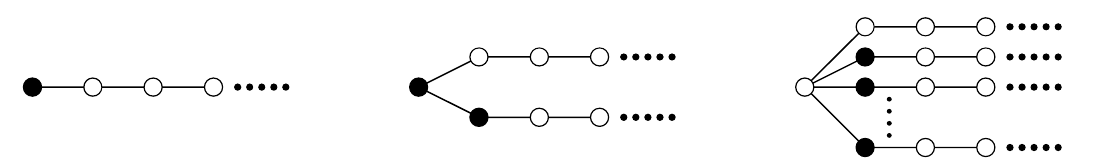}
\caption{Black vertices form a metric basis of the graphs
$P_\infty$ (left) , $P_{2\infty}$ (middle) and $P_{k\infty}$
(right). }\label{fig.infinitePaths}\end{center}
\end{figure}

A natural starting point for studying the metric dimension of
infinite graphs is to investigate the finiteness of this
invariant. It is not difficult to realize that infinite graphs may
have finite or infinite  metric dimension. In fact, for every
$k\ge 0$ there exist infinite graphs with metric
dimension $k$. To prove this, we define the \emph{one-way infinite
path}, $P_\infty $, as the infinite graph with set of vertices
$V(P_\infty)=\{ u_i : i \ge 0 \}$, and two vertices $u_i$, $u_j$
are adjacent if and only if $|i-j|=1$. We say that $u_0$ is the
\emph{endpoint} of $P_\infty$. In a similar way we define for
$k\ge 2$ the \emph{$k$-way infinite path}, $P_{k\infty}$, as the
graph formed by $k$ pairwise disjoint one-way infinite paths and a
new vertex adjacent to their $k$ endpoints. It is straightforward
to prove that $\beta (P_\infty )=1$, $\beta(P_{2\infty})=2$ and
$\beta( P_{k\infty})=k-1$, if $k\ge 3$ (see
Figure~\ref{fig.infinitePaths}). It is also easy to verify that
finite paths and the one-way infinite path  are the only
graphs with metric dimension equal to 1, since there exists a
vertex $v$ for which there is at most one vertex at distance $k$,
for every $k\ge 0$.


\begin{proposition}{\label{pro.dim1}}
Then the metric dimension of a graph $G$ is $1$ if, and
only if, $G$ is either a finite path 
or the one-way infinite
path.
\qed
\end{proposition}

\begin{figure}[ht]
\begin{center}
\includegraphics[width=0.28\textwidth]{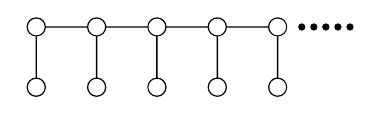}
\caption{The infinite comb graph $B_{\infty}$ has infinite metric
dimension.} \label{fig.infiniteComb}\end{center}
\end{figure}

So, the metric dimension of any infinite graph $G\not= P_{\infty
}$ is at least $2$. However, apart from $P_{2\infty}$ and
$P_{3\infty}$, there are many more infinite graphs  with metric
dimension $2$ as, for example, any graph obtained by attaching a
finite path or a one-way infinite path at an arbitrary  vertex of
degree $2$ of $P_{\infty}$. On the other hand, if we attach one
leaf at every vertex of the one-way infinite path  we obtain the
so-called \emph{infinite comb graph}, $B_{\infty}$,  which
satisfies $\beta(B_{\infty})=\infty$ (see
Figure~\ref{fig.infiniteComb}).


\begin{proposition}
The infinite comb graph $B_{\infty}$ has infinite metric
dimension.
\end{proposition}
\begin{proof}
Suppose that  $V(B_{\infty})=\{u_i : i\ge 0\} \cup \{ v_i : i\ge 0
\}$ and  $E(B_{\infty})=\{ u_iu_{i+1}: i\ge 0 \} \cup \{ u_iv_i:
i\ge 0 \}$. If $B_{\infty}$ has a finite resolving set $S$, there
exists a vertex $u_k$ such that $u_i$, $v_i$ are not in $S$ for
all $i\ge k$. Then, $d(w,u_{k+1})=d(w,u_k)+1=d(w,v_k)$ for every
$w\in S$. Hence,  $S$ does not resolve the pair $u_{k+1}$ and
$v_k$.
\end{proof}

In going on from finite to infinite graphs, a natural technique is
to study the desired property or parameter for their induced
subgraphs. Nevertheless, this seems to go nowhere in the case of
the metric dimension. Let us illustrate this remark with some
examples.

The infinite comb graph has infinite metric dimension,
nevertheless the metric dimension of every finite induced subgraph
is at most 2.

\begin{figure}[ht]
\begin{center}
\includegraphics[width=\textwidth]{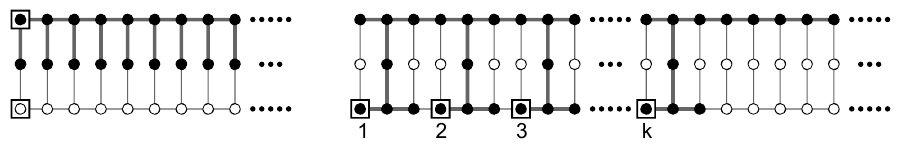}
\caption{Left, squared vertices form a metric basis of the illustrated
graph, while black vertices induce the graph $B_{\infty}$, with
infinite metric dimension. Right, top vertices induce the graph
$P_{\infty}$
 with metric
dimension 1, and for every $k\ge 2$, squared vertices form a
metric basis of the infinite subgraph induced by black vertices.} \label{fig.inducedSubgraph}\end{center}
\end{figure}

The metric dimension of the graph  in
Figure~\ref{fig.inducedSubgraph} is 2, nevertheless it does
contain induced subgraphs with infinite metric dimension, as for
example the infinite comb graph, and induced subgraphs with metric
dimension $k$, for every $k\ge 1$ (see
Figure~\ref{fig.inducedSubgraph}).

We conclude this section by obtaining two necessary conditions for
an infinite graph to have finite metric dimension.

\subsection{Uniformly locally finite graphs}\label{unlofi}

\begin{figure}[ht]
\begin{center}
\includegraphics[width=0.45\textwidth]{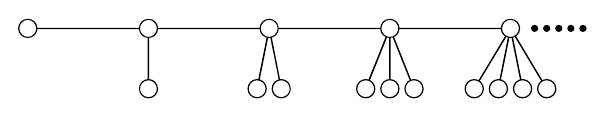}
\caption{An example of non uniformly locally finite graph.}
\label{fig.noULF}\end{center}
\end{figure}

An infinite graph is \emph{uniformly locally finite} if there
exists a positive integer $M$ such that the degree of every vertex
is at most $M$. For example, the graphs $P_{k\infty}$ and
$B_{\infty}$ are uniformly locally finite, and the graph obtained
by hanging $i$ vertices of degree $1$ to the vertex at distance
$i\ge 0$ from the endpoint of $P_{k\infty}$ is non uniformly
locally finite (see Figure~\ref{fig.noULF}).

\begin{lemma}\label{lem.degreeULF}
If $G$ is an infinite graph with
$\beta(G)=k$,
then every vertex of $G$ has degree at most $3^k-1$.
\end{lemma}
\begin{proof} Let $S=\{ x_1,\dots ,x_k\}$ be a metric basis of $G$.
Consider a vertex $v$ and its metric coordinates respect to $S$,
$r(v|S)=(d(v,x_1),d(v,x_2),\dots ,d(v,x_k))$. If $w$ is adjacent
to $v$, then $r(v|S)\not= r(w|S)$ and $|d(v,x_i)-d(w,x_i)|\le 1$
for all $x_i\in S$. This implies that there are only $3^k-1$
different possibilities for $r(w|S)$. Since all vertices must have
distinct metric coordinates, the degree of $v$ is at most $3^k-1$.
\end{proof}

As an immediate consequence of this lemma we get the following
necessary condition for an infinite graph to have finite metric
dimension.

\begin{theorem}\label{thm.uniflocally}
If $G$ is an infinite graph with finite metric dimension, then it
is uniformly locally finite.\qed
\end{theorem}

Note that the reciprocal of this  result is not true, being the
infinite comb graph illustrated in Figure~\ref{fig.infiniteComb} a
counterexample.

\subsection{Disjoint metric rays}\label{dimera}

A \emph{metric ray} with \emph{endpoint} $u_0$ in a graph $G$ is
an infinite subgraph $P$ of $G$ isomorphic to the one-way infinite
path, such that there exists an ordering of its vertex set,
$V(P)=\{ u_0,u_1,u_2,u_3,\dots \}$, with $u_k$ adjacent to
$u_{k+1}$ in $P$ for all $k\ge 0$, and
$d_G(u_0,u_k)=d_P(u_0,u_k)=k$.


Probably, the oldest and best known result about infinite graphs
is the K\"onig's infinity Lemma (see  \cite{konig}), that is
equivalent to the following assertion.

\begin{lemma}
For every vertex $v$ of an infinite graph,
there exists a metric ray with $v$ as endpoint.
\end{lemma}

Next, we show some results concerning metric rays which will help
us to find a second necessary condition for an infinite graph to
have finite metric dimension, quite different from that of
Theorem~\ref{thm.uniflocally}.

\begin{lemma}\label{lem.xk}
Let $P$ be a metric ray of an infinite graph $G$ with
ordered vertex set $V(P)=\{u_0,u_1,u_2,\ldots \}$.
Then, for every vertex  $x$ of $G$ there always exists an integer
$i_0\ge 0$ such that for every $j\ge i_0$,
$d(u_{j+1},x)=d(u_j,x)+1$.
\end{lemma}
\begin{proof}
Let $(a_i)_{i\ge 0}$ be the sequence defined as $a_i=
d(u_i,x)+d(x,u_0)-i$. Since $P$ is a metric ray, we have that
$d(u_i,x)+d(x,u_0)\ge d(u_i,u_0)=i$, hence $a_i\ge 0$. On the
other hand,  since $u_{i+1}$ and $u_i$ are adjacent vertices, we
have $d(u_{i+1},x)\le d(u_i,x)+1$, and consequently
$a_{i+1}=d(u_{i+1},x)+d(x,u_0)-(i+1)\le
d(u_{i},x)+1+d(x,u_0)-i-1=d(u_{i},x)+d(x,u_0)-i=a_i $. Hence,
$(a_i)_{i\ge0}$ is a decreasing sequence of non-negative integer
numbers. Therefore, there exists $i_0$ such that $a_{j+1}=a_j$ for
all $j\ge i_0$, or equivalently $a_j=k$ for all $j\ge i_0$. But
this implies that $d(u_{j+1},x)=d(u_j,x)+1$ for all $j\ge i_0$.
\end{proof}

\begin{lemma}\label{lem.uP}
Let $P$ be a metric ray of an infinite graph $G$, with
ordered vertex set $V(P)=\{u_0,u_1,u_2,\ldots \}$, and  $S$  a
finite subset of vertices of $G$.
There always exists an integer ${i_0}\ge 0$ such that for every
 $k\ge 0$,
$r(u_{i_0+k}|S)=r(u_{i_0}|S)+\underbrace{(k,\dots ,k)}_{|S|}$.

\end{lemma}
\begin{proof}
Let $S=\{x_1,\dots ,x_n\}$. By Lemma~\ref{lem.xk}, for every
$h\in[1,n]$ there exists ${i_h}$ such that
$d(u_{j+1},x_h)=d(u_j,x_h)+1$ for all $j\ge i_h$. Let $i_0=\max \{
{i_h} : h\in [1,n] \}$. Then $d(u_{j+1},x_h)=d(u_j,x_h)+1$ for all
$j\ge i_0$ and $h\in [1,n]$. This implies that
$d(u_{j+k},x_h)=d(u_j,x_h)+k$ for all $k\ge 0$, $j\ge i_0$ and
$h\in [1,n]$. In particular,
$r(u_{i_0+k}|S)=r(u_{i_0}|S)+\underbrace{(k,\dots ,k)}_{|S|}$.
\end{proof}

For every finite set $S$ of vertices and every metric ray $P$  in
an infinite graph $G$, we denote by ${i(P,S)}$ the
minimum among all  those integers $i_0$ satisfying the conditions
of the preceding lemma.

At this point, we recall the notion
 of doubly resolving set
introduced in~\cite{ours}. Two vertices $x$ and $y$ in a
nontrivial graph $G$  \emph{doubly resolve} a pair of
vertices $u$ and $v$ if $d(u,x)-d(v,x)\neq d(u,y)-d(v,y)$. If $S$
and $U$ are two subsets of vertices of $G$, we say that $S$
\emph{doubly resolves} $U$ if every pair of distinct vertices in
$U$ are doubly resolved by two vertices in $S$. We say that $S$ is
a \emph{doubly resolving set} of $G$ if $S$ doubly resolves
$V(G)$. If $G$ has at least one finite doubly resolving set, we
define $\psi (G)$ as the minimum cardinality of a doubly resolving
set; otherwise, we say that $\psi (G)=\infty $. Observe that every
doubly resolving set is also a resolving set, which means that
$\beta(G)\le \psi (G)$. Therefore, the next logical step should be
to study the finiteness of this parameter.


\begin{lemma}
A finite set of vertices does not doubly resolve any infinite set
of vertices of an infinite graph.
\end{lemma}

\begin{proof}
Suppose that $S=\{ x_1,...,x_n \}$ is a finite set and $U$ an
infinite set of vertices of a graph $G$. Let $D=\max \{ d(x_i,x_j)
: x_i,x_j\in S \}$. For every vertex $u\in U$, we define
$f_{(x_i,x_j)}(u)=d(u,x_i)-d(u,x_j)\in [-D,D]$,
 and $f_S(u)=(f_{p_1}(u),\dots
,f_{p_m}(u))$, where $\{ p_1,\dots ,p_{m}\}=\{(x_i,x_j):x_i,x_j\in
S, x_i\not= x_j\}$, $m=n(n-1)$. Since there are only $2D+1$
possible values for $f_{p_h}(u)=f_{(x_i,x_j)}(u)$, we have that
$|\{f_S(u): u\in U \} |\le (2D+1)^{n(n-1)}$. But $U$ in an
infinite set, so there are at least two different vertices $u,v\in
U$ such that $f_S(u)=f_S(v)$. This implies that for every pair
$(x_i,x_j)\in S\times S$, $x_i\not= x_j$, we have
$d(x_i,u)-d(x_j,u)=f_{(x_i,x_j)}(u)=f_{(x_i,x_j)}(v)=d(x_i,v)-d(x_j,v)$.
Then, $S$ does not doubly resolve $U$.
\end{proof}

\begin{corollary}\label{cor.psiInfinit}
If $G$ is an infinite graph, then $\psi (G)=\infty $.
\end{corollary}

\begin{lemma}
Let $G$ be an infinite  graph and let $S$ be a finite
resolving set of $G$. Suppose that $\mathcal{P}$ is a set of
metric rays in $G$ with pairwise disjoint vertex sets.
Then $S$ doubly resolves the set of vertices $\{ u_{i(P,S)} | P\in
\mathcal{P} \}$.
\end{lemma}
\begin{proof}
Suppose on the contrary that $S$  does not doubly resolve $W=\{
u_{i(P,S)} | P\in \mathcal{P} \}$. Then, there exists a pair of
distinct vertices $u,v\in W$ such that
$d(x,u)-d(x,v)=d(y,u)-d(y,v)$ for every pair of distint vertices
$x,y\in S$. This implies that $d(x,u)-d(x,v)=k$ for every $x\in
S$, and consequently $r(v|S)=r(u|S)+(k,\dots ,k)$. If $k\ge 0$, by
Lemma~\ref{lem.uP} there exists a vertex $w$ in the same metric
ray of $u$ such that $r(w|S)=r(u|S)+(k,\dots ,k)=r(v|S)$, which
contradicts the fact that $S$ is a resolving set of $G$. If $k<0$
we consider $d(x,v)-d(x,u)=-k>0$ and proceed in the same way.
\end{proof}

As a consequence of the preceding lemmas we obtain another
necessary condition for infinite graphs having finite metric
dimension.

\begin{theorem}\label{thm.charact}
If $G$ is an infinite  graph with finite metric
dimension, then it does not contain an infinite number of metric
rays with pairwise disjoint vertex sets.
\end{theorem}

\begin{figure}[ht]
\begin{center}
\includegraphics[width=0.28\textwidth]{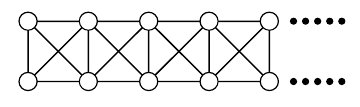}
\caption{A graph with infinite metric dimension not containing
infinite metric rays with pairwise disjoint vertex sets.}
\label{fig.finiteRays}\end{center}
\end{figure}

The reciprocal of Theorem~\ref{thm.charact} is not true: it is
easy to verify that the graph of Figure~\ref{fig.finiteRays} has
infinite metric dimension, however it contains at most two metric
rays with pairwise disjoint vertex sets.

\begin{figure}[ht]
\begin{center}
\includegraphics[width=0.25\textwidth]{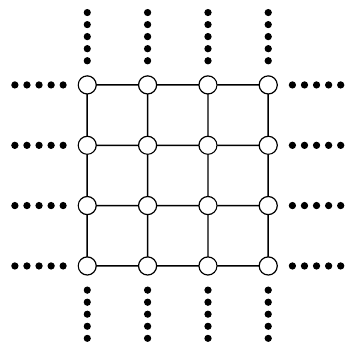}
\caption{The infinite grid.} \label{fig.infiniteGrid}\end{center}
\end{figure}

Note that the conditions stated in Theorems~\ref{thm.uniflocally}
and \ref{thm.charact} are independent. For example, the graph in
Figure~\ref{fig.noULF} is a non-uniformly locally finite graph not
containing infinite metric rays with pairwise disjoint vertex
sets. On the other hand, the \emph{infinite grid}, illustrated in
Figure~\ref{fig.infiniteGrid}, is uniformly locally finite, but
contains infinite metric rays with pairwise disjoint vertex sets.
Observe also that both conditions together are not sufficient to
assure finite metric dimension. For example, the infinite comb
graph has infinite metric dimension and satisfies both conditions.

\section{Infinite Trees}\label{sec.trees}

This section is devoted to study the metric dimension of
\emph{infinite trees}, that is, connected acyclic infinite graphs.
We obtain a finiteness characterization and calculate the exact
value of this parameter.


\begin{theorem}\label{thm.arbres}
An infinite tree has finite metric dimension if and only if the
set of vertices of degree at least three is finite.
\end{theorem}
\begin{proof} The only infinite trees without vertices of degree at least three are
$P_{\infty}$ and $P_{2\infty}$, which satisfy  $\beta (P_{\infty}
)=1$ and $\beta (P_{2\infty}) =2$. Now assume that $T$ is an
infinite tree  with some vertex of degree at least three.

($\Rightarrow$) Suppose first that $S$ is a finite resolving set
for an infinite tree $T$. Consider the finite set $S^*$ formed by
all vertices lying in a path with its endpoints in $S$. If there
are infinite vertices with degree at least three, there exists a
vertex $w\notin S^*$ such that $\deg (w)\ge 3$. Since $S^*$
induces a connected subgraph in $T$ and $T$ is acyclic, there is
only one vertex $z\in S^*$ such that $d(w,z)=\min \{ d(w,x) : x\in
S^*\}$. Let $w'$ be the vertex adjacent to $w$ that lies in the
unique $w-z$ path. Since $\deg (w)\ge 3$, there are at least two
vertices $u$, $v$ different from $w'$ and adjacent to $w$, such
that all paths to the vertices in $S\subseteq S^*$ pass through
$w$. Hence, $r(u|S)=r(w|S)+(1,\dots ,1)=r(v|S)$, contradicting the
hypotheses that $S$ is a resolving set.

($\Leftarrow$) Now suppose that the set $W=\{ x\in V : \deg (x)\ge
3 \}$ is finite. Consider the finite set $S$ formed by all
vertices lying in a path with endpoints in $W$ and all vertices
adjacent to some vertex of $W$. We claim that $S$ is a resolving
set for $T$. Since $S$ induces a connected subgraph in $T$ and
$T$ is acyclic, for every vertex $u$ not in $S$ there is only one
vertex $z(u)\in S$ such that $d(u,z(u))=\min \{ d(u,x) : x\in S
\}$. By definition of $S$, $z(u)\notin W$. This implies that $\deg
(z(u))=2$. Let $u$, $v$ be distinct vertices not in $S$ and
consider the vertices $z(u),z(v)\in S$. If $z(u)\not= z(v)$, then
there is only a path joining $u$ and $v$ that passes trough $z(u)$
and $z(v)$, and $z(u)$ or $z(v)$ resolves the pair $u$, $v$. If
$z(u)= z(v)$, then $\{ u,v,z(u) \}$ induces a path in $T$ with
$z(u)$ as an endpoint. Hence, $z(u)$ resolves $u$ and $v$.
\end{proof}

\begin{figure}[ht]
\begin{center}
\includegraphics[width=0.85\textwidth]{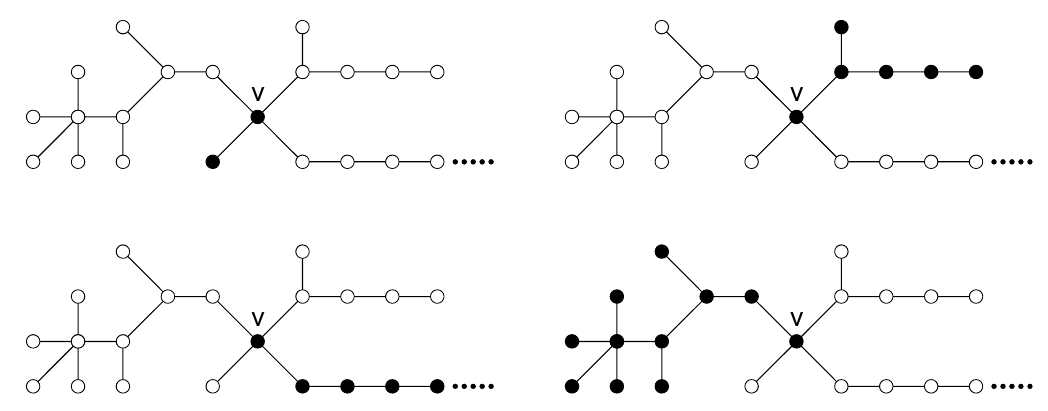}
\caption{Branches of an infinite tree $T$ at a vertex $v$ with
$\deg (v)=4$ and $P_T(v)=2$. At the left, the branch
paths.}\label{fig.branchPaths}\end{center}
\end{figure}

\begin{figure}[ht]
\begin{center}
\includegraphics[width=0.7\textwidth]{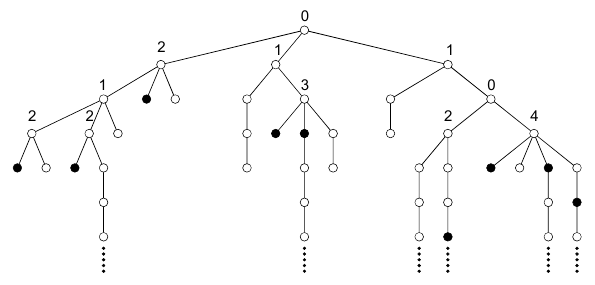}
\caption{The parameter $P_T(v)$ for vertices of degree at least
three. Black vertices form a metric basis of that infinite tree.
}\label{fig.baseArbreInf}\end{center}
\end{figure}

Let $v$ be a vertex of a (finite or infinite) tree $T$ with
maximum degree at least three. A \emph{branch} of $T$ at $v$ is a
maximal subtree having $v$ as a leaf. A \emph{branch path} of $T$
at $v$ is a branch that is either a finite path or a $1$-way
infinite path. Let $P_T(v)$ be the number of branch paths of $T$
at $v$ (see Figure~\ref{fig.branchPaths}). Observe that the number
of branches of $T$ at $v$ is exactly the degree of $v$, and
$P_T(v)\ge 2$ implies $\deg(v)\ge 3$. In \cite{chartrand,
landmarks, slater} it was proved that the metric dimension of a
finite tree different from a path is $\beta(T)=\sum_{\deg (v)\ge
3} \max \{ P_T(v)-1,0\}$. It was also shown that every metric
basis of $T$ is obtained by taking, for each vertex $v$ such that
$P_T(v)\ge 2$, exactly $P_T(v)-1$ vertices different from $v$
lying on different branch paths of $T$ at $v$. Next we prove that
the same results also hold for infinite trees (see
Figure~\ref{fig.baseArbreInf}).

\begin{theorem}
If  $T$ is an infinite tree with maximum degree at least three and
finite metric dimension, then
$$\beta(T)=\sum_{\deg (v)\ge 3} \max \{ P_T(v)-1,0\},$$
and all metric basis can be obtained in the following way: for
each vertex $v$ such that $P_T(v)=k\ge 2$, select $k-1$ vertices
different from $v$ lying on distinct branch paths of $T$ at $v$.
\end{theorem}

\begin{proof}
Let $S$ be one of the sets described in the theorem.  By
Theorem~\ref{thm.arbres}, the set $W=\{ v\in V(T) : \deg (v)\ge 3
\}$ is finite. Since $P_T(v)\ge 2$ implies $\deg (v)\ge 3$, the
set $S$ is also finite.

A resolving set of $T$ has at least so many vertices as $S$, since
for any $v\in W$ such that $P_T(v)\ge 2$, two vertices adjacent to
$v$ of different branch paths of $T$ at $v$ are resolved only by
vertices of those two branch paths and different from $v$. Observe
also that the sets of vertices of two different branch paths of
$T$ at vertices of degree at least three are disjoint. Therefore,
$|S|=\sum_{\deg (v)\ge 3} \max \{ P_T(v)-1,0\}$.

We claim that $S$ resolves $T$. Consider two distinct vertices
$x,y\in V(T)$. Let $T'$ be the graph induced by the set of
vertices that lie in some path with endpoints in $W\cup S\cup \{
x,y \}\cup N(W) $, where $N(W)$ is the set of  vertices adjacent
to some vertex of $W$. Then, $T'$ is a finite tree that contains
$W\cup S$,  for all $v\in W$ we have $P_T(v)=P_{T'}(v)$, and each
branch path of $T'$ at $v$  is contained in a branch path of $T$
at $v$. Therefore, $S$ is a resolving set for $T'$, and,
consequently, $S$ resolves the pair $x$, $y$ in $T$.
\end{proof}

\section{Cartesian
Products}\label{sec.cartProd}

\subsection{Conditions for finite metric dimension}

The cartesian product of graphs G and H, denoted by $G\Box H$, is
the graph with vertex set $V(G)\times V(H)=\{(a, v):a\in V(G),
v\in V(H)\}$, where $(a, v)$ is adjacent to $(b,w)$ whenever $a=b$
and ${v,w}\in E(H)$, or $v=w$ and ${a, b}\in E(G)$. Observe that
if $G$ and $H$ are connected, then $G\Box H$ is connected. In
particular, $d_{G\Box H}((a,v),(b,w))=d_G(a,b)+d_H(v,w)$. A number
of interesting examples can be viewed as the cartesian product of
two graphs. For example, the \emph{two dimensional infinite grid}
is the graph $P_{2\infty}\Box P_{2\infty}$.

Since the cartesian product of two infinite graphs has always
infinite metric rays with pairwise disjoint vertex sets, it
follows from Theorem~\ref{thm.charact} that:

\begin{corollary}
Given two infinite graphs $G$ and $H$, then $\beta(G\Box
H)=\infty$.
\end{corollary}

Given a subset $S$ of vertices in $G\Box H$, its \emph{projection}
onto $G$ is the set of vertices $x \in V(G)$ for which there
exists a vertex $u\in V(H)$ such that $(x,u)\in S$. Similarly is
defined the projection of $S$ onto $H$. Some results obtained for
finite graphs  in ~\cite{ours} can be extended to infinite graphs.

\begin{proposition}
\label{prop.proj} Let $G$, $H$ be finite or infinite graphs.
If $S$ is a resolving set of $G\Box H$, then the projection of $S$
onto $G$ (resp. onto $H$) is a resolving set of $G$ (resp. of $H$)
and, consequently,  $\beta (G\Box H)\ge \max \{ \beta (G), \beta
(H)\}$.
\end{proposition}

\begin{corollary}
If $G$ is an infinite graph with infinite metric
dimension, then for any graph $H$ we have $\beta(G\Box H)=\infty$.
\end{corollary}

Another result that can be extended to infinite graphs is the
following:

\begin{theorem}
\label{thm.psi}  If $G$ is an infinite graph with finite metric dimension and $H$ is a finite
graph with at least two vertices, then the metric
dimension of $G\Box H$ is finite and $\beta(G\Box H)\le \beta (G)
+\psi (H)-1$
\end{theorem}

We do not include the proofs of Proposition \ref{prop.proj} and
Theorem \ref{thm.psi}, since the proofs given in \cite{ours} also
hold for infinite vertex sets. \vspace{0.3cm}

We summarize the preceding results in Table~\ref{tab.cartProduct}.


\begin{table}[h]
 \begin{center}
 \begin{tabular}{|cc|c|}
  \hline
  $\phantom{x}|V(G)|\phantom{x}$ & $\phantom{x}|V(H)|\phantom{x}$ & $\phantom{x}\beta (G\Box H)\phantom{x}$ \\
  \hline
  $<\infty $ & $<\infty $ & $<\infty $\\
  \hline
  $=\infty $ & $=\infty $ & $=\infty $ \\
  \hline
  $=\infty $ & $<\infty $ &
  $\begin{cases}
    <\infty  & \text{ if }\beta (G)<\infty,\\
    =\infty  & \text{ if }\beta (G)=\infty.
  \end{cases}$\\
 \hline
 \end{tabular}
 \end{center}
 \caption{Possibilities for the metric dimension of the cartesian product of graphs.}
 \label{tab.cartProduct}       
\end{table}

\subsection{Metric dimension of $P_{\infty}\Box H$ and  $P_{2\infty}\Box H$}

In this section, we study the metric dimension of the cartesian
product of the infinite graphs $P_{\infty}$ and $P_{2\infty}$ by
finite graphs. In particular we determine the metric dimension of
$P_{\infty}\Box H$ and $P_{2\infty} \Box H$ when $H$ is a path, a
cycle  or a complete graph. Proposition ~\ref{prop.proj} gives us
a lower bound of $\beta (G\Box H)$ in terms of $\beta (H)$. It is
known that the metric dimension of paths, cycles and complete
graphs of order $n\ge 3$ is respectively $1$, $2$ and $n-1$ (see
~\cite{chartrand,landmarks}). Moreover, a vertex of degree $1$ is
a basis of a path; any two no antipodal vertices of a cycle form a
basis of a cycle, and a basis of a complete graph is formed by any
$n-1$ vertices of the graph. On the other hand, Theorem
~\ref{thm.psi} gives an upper bound of $\beta (G\Box H)$ in terms
of $\psi (H)$. This parameter is determined in \cite{ours} for
several families of graphs as paths $P_n$, cycles $C_n$ and
complete graphs $K_n$ of order $n$. Concretely, $\psi (P_n)=2$, if
$n\ge 2$; $\psi (C_{n})=2$, if $n\ge 3$ odd; $\psi(C_{n})=3$, if
$n\ge 4$ even; and $\psi (K_n)=n-1$, if $n\ge 3$. In the rest of
the paper we suppose that $V(P_n)=V(C_n)=V(K_n)=\{0,1,2,\dots ,
n-1\}$, $V(P_{\infty}) =\mathbb{N}$ and $V(P_{2 \infty})
=\mathbb{Z}$. Two vertices $i$, $j$ are adjacent in $P_n$ or
$P_\infty $ or $P_{2\infty }$ if and only if $|j-i|=1$.  Two
vertices $0\le i\le j\le n-1$ are adjacent in $C_n$ if and only if
$j-i=1$ or $j-i=n-1$, and any two distinct vertices of $K_n$ are
adjacent.

\begin{table}[h]
 \begin{center}
 \begin{tabular}{|c|cccc|}
  \hline
   $G$& \phantom{m}$P_n$ & \phantom{m}$C_{n}$, $n$ odd\phantom{m} &
    \phantom{m}$C_{n}$, $n$ even\phantom{m} & \phantom{m}$K_n$\phantom{m}   \\
  \hline
 $\beta(G)$& 1&2&2& $n-1$ \\
 $\psi(G)$& 2&2&3& $n-1$\\
  \hline
 \end{tabular}
 \end{center}
 \caption{Values of  $\beta(G)$ and $\psi(G)$ of paths, cycles and cliques of order $n\ge3$.}
 \label{tab.finitecase}
\end{table}

\begin{lemma}\label{lem.p2infinitH}
If $H$ is a  graph, then $\beta(P_{2\infty}\Box H)=2$ if
and only if $H$ is the trivial graph, $K_1$.
\end{lemma}
\begin{proof}
($\Leftarrow$) The graph $P_{2 \infty} \Box K_1$ is  $P_{2
\infty}$ with metric dimension  $2$.

($\Rightarrow$) Suppose now that $H$ is a non trivial
graph $H$ and  $S=\{ (i,u),(j,v)\}$, $i\le j$, is a  set with two
different vertices of  $P_{2\infty}\Box H$. We claim that $S$ does
not resolve $P_{2\infty}\Box H$. By Proposition ~\ref{prop.proj},
if $i=j$, then the projection of $S$ onto $P_{2\infty}$ is not a
resolving set of $P_{2\infty}$. Now suppose that $i\not= j$.

Case  $u=v$. Let $w$ be a vertex adjacent to $u$ in $H$, that
exists because $H$ is a non trivial graph. Then
$d_{P_{2\infty}\Box H}((i-1,u),(i,u))=1=d_H(w,u)=d_{P_{2\infty}
\Box H }((i,w),(i,u))$ and $d_{P_{2\infty}\Box
H}((i-1,u),(j,u))=d_{P_{2\infty}}(i-1,j)=d_{P_{2\infty}}(i,j)+1=d_{P_{2\infty}}(i,j)+d_H(w,u)=d_{P_{2\infty}\Box
H}((i,w),(j,u))$. Hence, $S$ does not resolve the vertices
$(i-1,u)$ and $(i,w)$.

Case $u\not= v$. Let $w$ be a vertex of $H$ adjacent to $u$ and
lying in a shortest $u-v$ path. Then
 $d_{P_{2\infty}\Box
H}((i+1,u),(i,u))=1=d_H(w,u)=d_{P_{2\infty} \Box H}((i,w),(i,u))$
and $d_{P_{2\infty}\Box
H}((i+1,u),(j,v))=d_{P_{2\infty}}(i+1,j)+d_H(u,v)=d_{P_{2\infty}}(i,j)-1+d_H(u,v)=
d_{P_{2\infty}}(i,j)+d_H(w,v)=d_{P_{2\infty}\Box H}((i,w),(j,v))$.
Hence, $S$ does not resolve the vertices $(i+1,u)$ and $(i,w)$.
\end{proof}


\begin{lemma}\label{lem.Smesv}
If $H$ is a graph and $S\subseteq \{0 \}\times V(H)$ is
a resolving set of $P_{\infty}\Box H$, then, for any $u\in V(H)$,
$S'=S\cup \{ (1,u)\}$ is a resolving set of $P_{2\infty}\Box H$.
\end{lemma}
\begin{proof}
Let $x=(i,v)$ and $y=(j,w)$ be two distinct vertices of
$P_{2\infty }\Box H$. If $i,j\ge 0$, $S$ resolves $x$ and $y$. By
symmetry, $S$ resolves $x$ and $y$ if $i,j\le 0$.

Suppose now that $i>0$ and $j<0$. Since $r((j,w)|S)=r((-j,w)|S)$
and $S$ is a resolving set for $P_{\infty}\Box H$, $S$ resolves
$x$ and $y$ if $i\not= -j$ or $v\not= w$. Finally, for $i=-j>0$
and $v=w$ we have $d_{P_{2\infty }\Box H}(x,(1,u))=d_{P_{2\infty
}\Box H}((i,v),(1,u))=i-1+d_H(v,u)\not= i+1+d_H(v,u)
=d_{P_{2\infty }\Box H}((-i,v),(1,u))= d_{P_{2\infty }\Box
H}((j,w),(1,u))= d_{P_{2\infty }\Box H}(y,(1,u))$. Hence, $(1,u)$
resolves $x$ and $y$, and consequently $S'$ is a resolving set for
$P_{2\infty}\Box H$.
\end{proof}


We determine now the metric dimension and a metric basis of $G\Box
P_n$, when $G$ is $P_{\infty}$ or $P_{2\infty }$ and $n\ge 2$. We
have in that case $$d_{G\Box
P_n}((i,j),(i',j'))=d_G(i,j)+d_{P_n}(i',j')=| i'-i |+|j'-j|.$$

\begin{proposition}\label{pro.pinfinitcami} For all $n\ge 2$, $\beta(P_\infty \Box P_n)=2$
and  $S=\{ (0,0), (0,n-1) \}$ is a metric basis of $P_\infty \Box
P_n$.
\end{proposition}
\begin{proof}
By Theorem ~\ref{thm.psi}, $\beta(P_\infty \Box P_n)\le
\beta(P_\infty)+\psi (P_n)-1=1+2-1=2$ if $n\ge 2$. Since the
infinite graph $P_\infty \Box P_n$ is not the one-way infinite
path for $n\ge 2$, by Proposition \ref{pro.dim1} we have $
\beta(P_\infty \Box P_n)=2$ for all $n\ge 2$.

Consider now two different vertices $(i,j)$ and $(i',j')$ in
$P_{\infty}\Box P_n$. $S$ is a resolving set of the subgraph
isomorphic to $P_m\Box P_n$ that contains $S$, $(i,j)$ and
$(i',j')$, where $m=\max \{ i,i'\}+1$ (see \cite{landmarks}).
Then, $S$ resolves the pair $(i,j)$ and $(i',j')$ in
$P_{\infty}\Box P_n$.
\end{proof}

\begin{figure}[h]
\begin{center}
\includegraphics[width=0.5\textwidth]{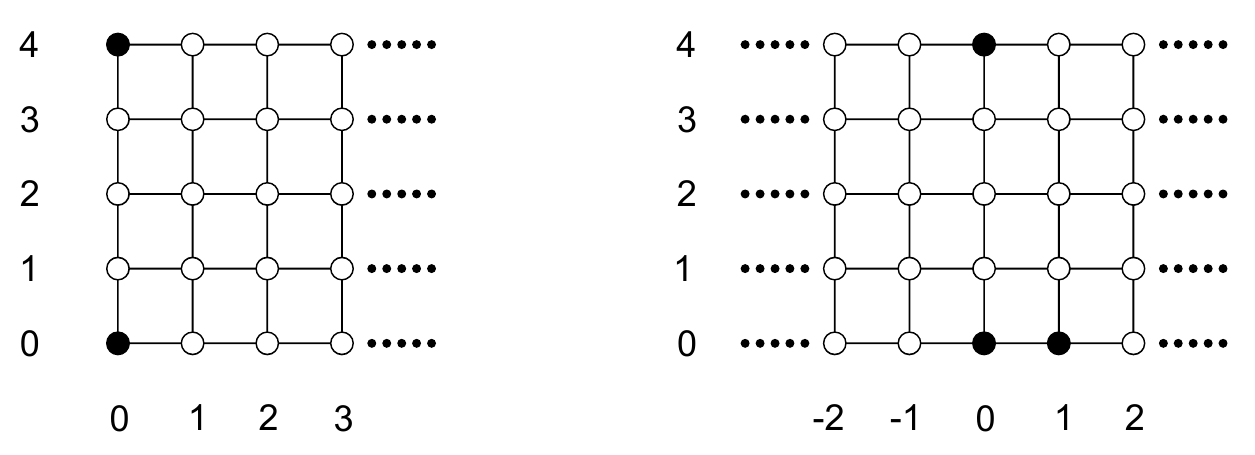}
\caption{Black vertices form a metric basis of the graphs
$P_\infty \Box P_5$ (left) and $P_{2\infty} \Box P_5$ (right).
}\label{fig.baseCamins}
\includegraphics[width=0.4\textwidth]{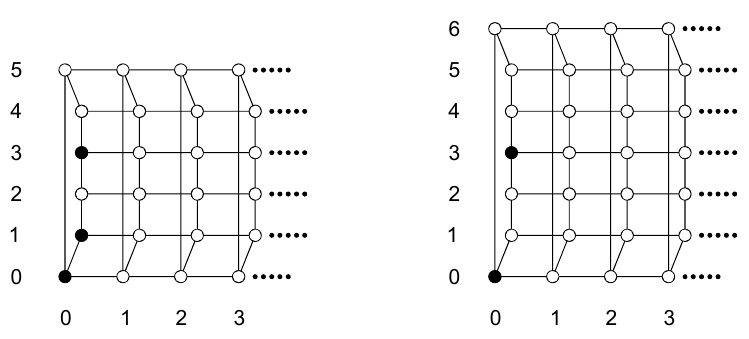}\hspace{0.2cm}
\includegraphics[width=0.5\textwidth]{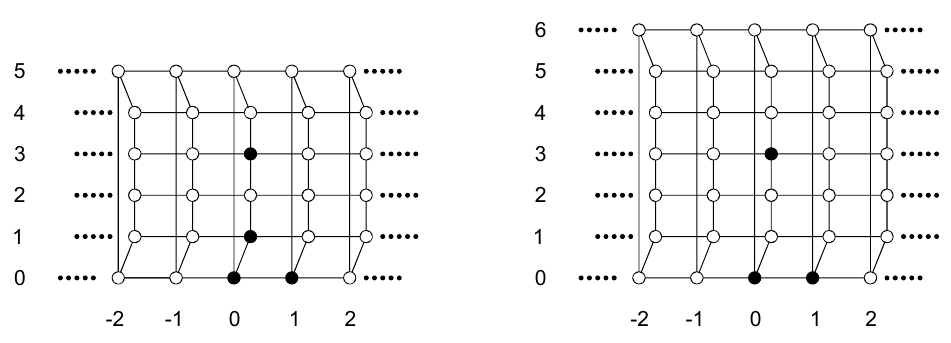}
\caption{From left to right, black vertices form a metric basis of the graphs
$P_\infty \Box C_6$
, $P_{\infty} \Box C_7$
, $P_{2\infty} \Box C_6$
and $P_{2\infty} \Box C_7$.}
\label{fig.baseCicles}
\vspace{0.5truecm}
\includegraphics[width=0.5\textwidth]{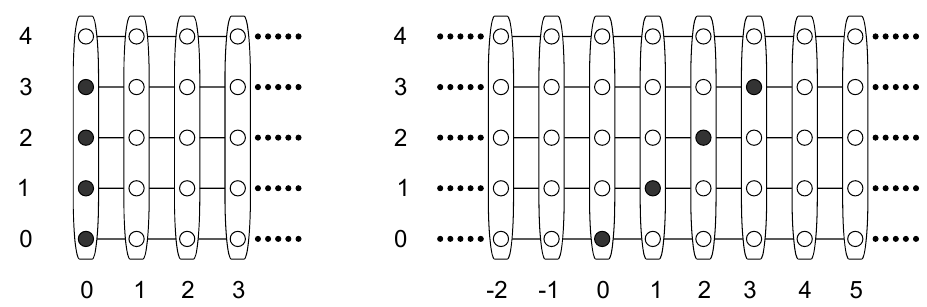}
\caption{Black vertices form a metric basis of the graphs
$P_\infty \Box K_5$ (left) and $P_{2\infty} \Box K_5$ (right). The
vertices of each column are pairwise
adjacent.}\label{fig.baseComplets}
\end{center}
\end{figure}

\begin{proposition}\label{pro.p2infinitcami} For all $n\ge 2$,
$\beta(P_{2\infty}\Box P_n)=3$  and $S=\{ (0,0), (0,n-1),(1,0)\}$
is a metric basis of  $P_{2\infty} \Box P_n$.
\end{proposition}
\begin{proof}
By Theorem~\ref{thm.psi} and Proposition ~\ref{prop.proj},
$2=\beta(P_{2\infty })\le \beta(P_{2\infty}\Box P_n)\le \beta(P_{2
\infty})+\psi (P_n)-1=2+2-1=3$. By Lemma ~\ref{lem.p2infinitH},
$\beta(P_{2\infty}\Box P_n)\not= 2$. Hence, $\beta(P_{2\infty}\Box
P_n)=3$. See Figure~\ref{fig.baseCamins}.

By  Proposition~\ref{pro.pinfinitcami} and Lemma~\ref{lem.Smesv},
$S$ is a metric basis of $P_{2\infty}\Box P_n$.
\end{proof}

For determining the metric dimension and a metric basis of $G\Box
C_n$, when $G$ is $P_{\infty}$ or $P_{2\infty }$ and $n\ge 3$,
observe that
$$d_{G\Box C_n}((i,j),(i',j'))=d_G(i,j)+d_{C_n}(i',j')=| i'-i
|+\min \{ |j'-j|,n-|j'-j|\}.$$

\begin{proposition}\label{pro.pinfinitCicle} For all $n\ge 3$, the
metric dimension of  $P_\infty \Box C_n$ is
$$\beta(P_\infty \Box C_n)=
  \begin{cases}
    2 & \text{if $n$ is odd}, \\
    3 & \text{if $n$ is even},
  \end{cases}$$
  and $S_1=\{ (0,0), (0,\frac{n-1}2)\}$ is a metric basis, if $n$ is odd,
 and $S_2=\{ (0,0), (0,\frac n2),(0,1)\}$ is a metric basis,  if $n$ is even.
\end{proposition}
\begin{proof}
By Theorem~\ref{thm.psi} and Proposition~\ref{prop.proj}, $2=\beta
(C_n)\le  \beta(P_\infty \Box C_n)\le \beta(P_{\infty})+\psi
(C_n)-1=\psi (C_n)$.
Therefore, $\beta(P_\infty \Box C_n)$ is $2$, if $n$ is odd, and
$2$ or $3$, if $n$ is even.

Suppose now that $n$ is even and $W=\{ (i,j), (i',j')\}$
 is a set with two different vertices of $P_\infty \Box C_n$. If
$W$ is a resolving set of $P_\infty \Box C_n$, the projection of
$W$ onto $P_\infty$ is a resolving set of $P_\infty$ and the
projection onto $C_n$ is a resolving set of $C_n$. Therefore,
$i=i'=0$ or $i\not= i'$, and $j$, $j'$ are different vertices not
antipodal in $C_n$. By symmetry, we may assume that $0=j< j'<
\frac n2$. If $i=i'=0$, $W$ does not resolve the pair of vertices
$(0,j'+1)$ and $(1,j')$ (see Figure ~\ref{fig.prop6}, left). If
$i\not= i'$, we may assume by symmetry that $0\le i<i'$. Then $W$
does not resolve the pair of vertices $(i,1)$ and $(i+1,0)$  (see
Figure ~\ref{fig.prop6}, right). In any case, we have a
contradiction. Hence $\beta (P_{\infty}\Box C_n)=3$ for $n$ even.

\begin{figure}[ht]
\begin{center}
\includegraphics[width=0.7\textwidth]{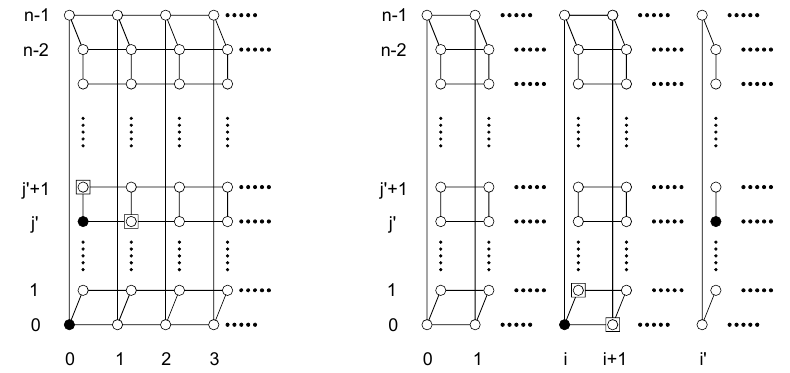}
\caption{For $n$ even, black vertices do not resolve squared
vertices if  $0<j'<\frac n2$ (left) and  if $0\le i<i'$,
$0<j'<\frac n2$ (right).}\label{fig.prop6}
\end{center}
\end{figure}

Now, for $n=2k+1$ odd,  consider the set $S_1=\{ (0,0), (0,k)\}$.
Then,

\begin{equation}\label{eqn.coordSenar}
r((0,j)|S_1)=
  \begin{cases}
    (j,k-j) & \text{ if }0\le j\le k, \\
    (2k+1-j,j-k) & \text{ if }k+1\le j\le 2k.
  \end{cases}
  \end{equation}

Observe that the sum of the two coordinates is $k$ in the first
case and $k+1$ otherwise. For $i> 0$ we have $r((i,j)|S_1)=
r((0,j)|S_1)+(i,i)$. Now, the sum of the two coordinates is $2i+k$
if $0\le j\le k$ and $2i+k+1$ if $k+1\le j\le 2k$. Suppose that
the vertices $x=(i,j)$ and $y=(i',j')$, where $0\le i\le i'$,
satisfy $r((i,j)|S_1)= r((i',j')|S_1)$. This implies
$r((0,j)|S_1)=r((i'-i,j')|S_1)$. But this is not possible when
$i'-i\not= 0$, since the sum of the two coordinates of
$r((0,j)|S_1)$ is $k$ or $k+1$, and the sum of the two coordinates
of $r((i'-i,j')|S_1)$ is $2(i'-i)+k$ or $2(i'-i)+k+1$. Then
 $i'=i$ and, from \eqref{eqn.coordSenar}, we have $j=j'$. Thus, $S_1$ resolves
$P_{\infty}\Box C_n$ if $n$ is odd.

If  $n=2k$ even,  consider the set $S_2=\{ (0,0), (0,k), (0,1)\}$.
Then,

\begin{equation}\label{eqn.coordParell}
r((0,j)|S_2)=
  \begin{cases}
   (0,k,1) & \text{ if $j=0$},\\
    (j,k-j,j-1) & \text{ if $1\le j\le k$}, \\
    (2k-j,j-k,2k-j+1) & \text{ if $k+1\le j\le 2k-1$}.\\
  \end{cases}
\end{equation}

Observe that the sum of the first two coordinates is $k$.  For $i>
0$ we have $r((i,j)|S_2)= r((0,j)|S_2)+(i,i,i)$, and the sum of
the first two coordinates is $2i+k$.

Suppose that the vertices $x=(i,j)$ and $y=(i',j')$, where $0\le
i\le i'$, satisfy $r((i,j)|S_2)= r((i',j')|S_2)$. This implies
$r((0,j)|S_2)=r((i'-i,j')|S_2)$. But this is not possible when
$i'-i\not= 0$, since the sum of the first two coordinates of
$r((0,j)|S_2)$ is $k$, and the sum of the first two coordinates of
$r((i'-i,j')|S_2)$ is $2(i'-i)+k$. Then  $i'=i$ and, from
\eqref{eqn.coordParell}, we have that $j=j'$. Thus, $S_2$ resolves
$P_{\infty}\Box C_n$ if $n$ is even.
\end{proof}

\begin{proposition}\label{pro.p2infinitCicle} For all $n\ge 3$,
$$\beta(P_{2\infty} \Box C_n)=
  \begin{cases}
    3 & \text{if $n$ is odd}, \\
    4 & \text{if $n$ is even}.
  \end{cases}$$
and a metric basis  of $P_{2\infty} \Box C_n$ is $S_1=\{ (0,0),
(0,\frac{n-1}2), (1,0)\}$, if $n$ is odd, and $S_2=\{
   (0,0), (0,\frac n2),(0,1),(1,0)\}$, if $n$ is even.
\end{proposition}
\begin{proof}
By Theorem~\ref{thm.psi} and Proposition~\ref{prop.proj},
$2=\beta(C_n)\le \beta(P_{2\infty }\Box C_n)\le
\beta(P_{2\infty}+\psi (C_n)-1=\psi (C_n)+1$.
By Lemma~\ref{lem.p2infinitH}, $\beta(P_\infty \Box C_n)\not=2$.
Therefore,  $\beta(P_{2\infty }\Box C_n)=3$, if $n$ is odd, and
$\beta (P_\infty \Box C_n)$  is $3$ or $4$, if $n$ is even. On the
other hand, Proposition~\ref{pro.pinfinitCicle} and
Lemma~\ref{lem.Smesv} imply that $S_1$ is a metric basis of
$P_\infty \Box C_n$ for $n$ odd.

\begin{figure}[ht]
\begin{center}
\includegraphics[width=1\textwidth]{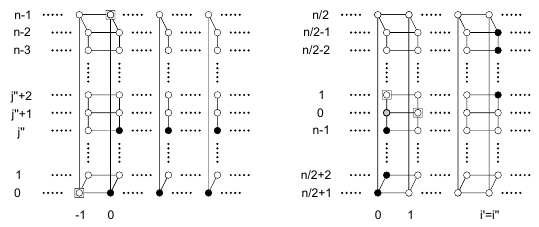}
\caption{ For $n$ even, three black vertices of different columns
do not resolve squared vertices if $0<j''<\frac n2$ (left) and
vertex $(0,0)$ together with  two black vertices of different
columns do not resolve squared vertices
(right).}\label{fig.prop7cas2j}
\end{center}
\end{figure}
\begin{figure}[ht]
\begin{center}
\includegraphics[width=\textwidth]{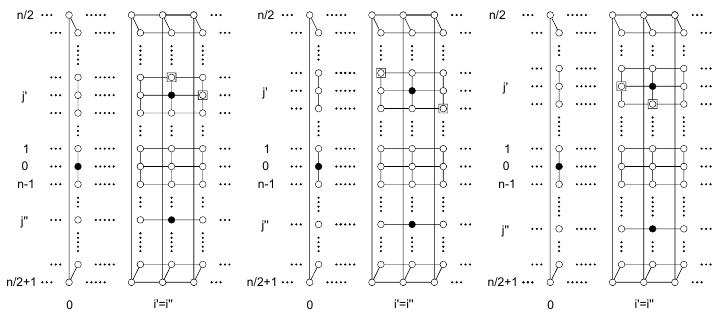}
\caption{Black vertices do not resolve squared vertices if
$j''-j'>n/2$ (left), if $j''-j'=n/2$ (middle), and if $j''-j'<n/2$
(right).}\label{fig.prop7cas3j2ib}
\end{center}
\end{figure}

Suppose now that $n$ is even and let $S=\{ (i,j),(i',j'),
(i'',j'')\}$ be a set vertices of $G=P_{2\infty }\Box C_n$ of
cardinality $3$. We claim that $S$ does not resolve $G$.

If $|\{ j,j',j''\} |=1$ or  $\{ j,j',j''\}$ is a set of two
antipodal vertices of $C_n$, then $S$ does not resolve $P_\infty
\Box C_n$ by Proposition~\ref{prop.proj}.

If $\{ j,j',j''\}$ is a set of two no antipodal vertices of $C_n$,
we may assume by symmetry that $j=j'=0$, $0<j''<\frac n2$ and
$\min \{ i,i',i''\}=0$. Then $S$ does not resolve $(-1,0)$ and
$(0,n-1)$ (see Figure ~\ref{fig.prop7cas2j}, left).

If $| \{ j,j',j'' \} | =3$, at least one vertex is not the
antipodal of the remaining two vertices. In any case, we may
assume without lose of generalization that $j=0$, $0<j'<n/2$ and
$n/2 < j'' \le n-1 $.
Now, if $|\{ i,i',i''\}|=1$,  $S$ is not a resolving set by
Proposition~\ref{prop.proj}.
If $|\{ i,i',i''\}|=2$, we may assume by symmetry that all the
cases are analogous to (1) $i=i''=0$ and $i'>0$, or (2) $i=0$ and
$i'=i''>0$. In the first case, $S$ does not resolve $(0,1)$ and
$(1,0)$ (see Figure~\ref{fig.prop7cas2j}, right).
In the second
case, $S$ does not resolve $(i',j'+1)$ and $(i'+1,j')$ if
$j''-j'>\frac n2$; $S$ does not resolve $(i'-1,j'+1)$ and
$(i'+1,j'-1)$ if $j''-j'=\frac n2$; and $S$ does not resolve
$(i'-1,j')$ and $(i',j'-1)$ if $j''-j'<\frac n2$ (see Figure
~\ref{fig.prop7cas3j2ib}).

Finally, if $|\{ i,i',i''\}|=3$, by symmetry all the cases are
analogous to (1) $i''<i=0<i'$ or (2) $i<i''=0<i'$. In both cases
$S$ does not resolve $(0,1)$ and $(1,0)$ (see Figure
~\ref{fig.prop7cas3j3i}).

\begin{figure}[ht]
\begin{center}
\includegraphics[width=\textwidth]{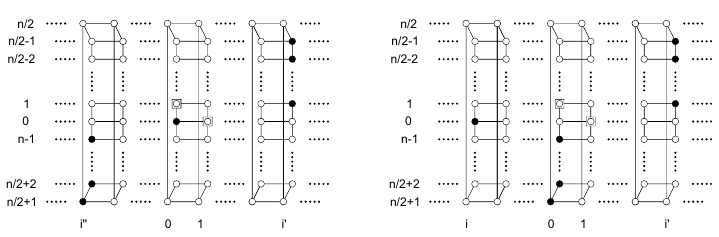}
\caption{Three black vertices of different columns do not resolve
squared vertices.}\label{fig.prop7cas3j3i}
\end{center}
\end{figure}

Thus, $S$ does not resolve $P_{2\infty }\Box C_n$ when $n$ is
even, implying that the metric dimension is 4. By
Proposition~\ref{pro.pinfinitCicle} and Lemma~\ref{lem.Smesv},
$S_2$ is a metric basis of $P_{2\infty}\Box C_n$.
\end{proof}


Finally, we give the metric dimension and a metric basis of $G\Box
K_n$, when $G$ is $P_{\infty}$ or $P_{2\infty }$ and $n\ge 4$. We
have in that case $$d_{G\Box K_n}((i,j),(i',j'))=
  \begin{cases}
     | i'-i |\, , & \text{ if } j=j', \\
     | i'-i |+1\, , & \text{ if } j\not= j'.
  \end{cases}
$$

\begin{proposition}\label{pro.pinfinitComplet} For all $n\ge 4$, $\beta(P_\infty \Box K_n)=n-1$
and $S=\{ (0,0),(0,1),\dots ,(0,n-2)\}$ is a metric basis of
$P_\infty \Box K_n$.
\end{proposition}
\begin{proof}

By Theorem~\ref{thm.psi} and Proposition~\ref{prop.proj},
$n-1=\beta(K_n)\le \beta(P_{\infty }\Box K_n)\le
\beta(P_{\infty})+\psi (K_n)-1=\psi (K_n)=n-1$. Hence,
$\beta(P_{\infty }\Box K_n)=n-1$.

Consider now a pair of two distinct vertices, $x=(i,j)$,
$y=(i',j')$. If $i=i'$, the vertex $(0,h)$ resolves $x$ and $y$,
where $h=\min \{ j,j'\}\in [0,n-2]$. If $i\not= i'$, the vertex
$(0,h)$ resolves $x$ and $y$, where $h$ is any integer $h\in
[0,n-2]$ such that $h\not= j,j'$. Thus, $S$ is a metric basis of
$P_{\infty }\Box K_n$.
\end{proof}

\begin{proposition}\label{pro.p2infinitComplet} For all $n\ge 4$, $\beta(P_{2\infty}\Box K_n)=n-1$
and $S=\{ (0,0),(1,1),\dots ,(n-2,n-2)\}$ is a metric basis of
$P_{2\infty} \Box K_n$.
\end{proposition}
\begin{proof}

By Theorem~\ref{thm.psi} and Proposition~\ref{prop.proj},
$n-1=\beta (K_n)\le \beta(P_{2\infty }\Box K_n)\le
\beta(P_{2\infty})+\psi (K_n)-1=\psi (K_n)+1=n$. Hence,
$\beta(P_{\infty }\Box K_n)$ is $n-1$ or $n$.

We show next that $S$ is a metric basis of $P_{2\infty }\Box K_n$.
Let $x=(i,j)$, $y=(i',j')$ be
two different vertices of $P_{2\infty}\Box K_n$. If $i=i'$, the
vertex $(h,h)$ resolves $x$ and $y$, where $h=\min \{ j,j'\}\in
[0,n-2]$. Suppose now that $i\not= i'$. Consider $k\in [0,n-2]$
such that $k\not=j,j'$ and with the additional condition $k\not=
\frac {i+i'}2$ if the distance $d_{P_{2\infty}}(i,i')$ is even.

\begin{figure}[ht]
\begin{center}
\includegraphics[width=0.6\textwidth]{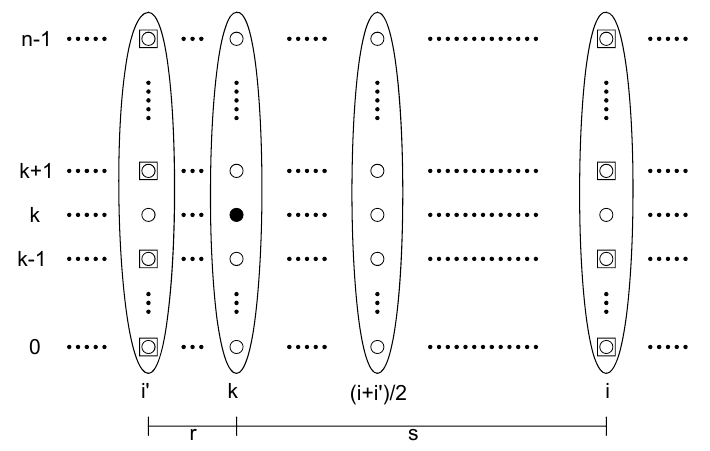}
\caption{Black vertex resolves two squared vertices  of different
columns if $r\not= s$.}\label{fig.prop9cas2i}

\includegraphics[width=\textwidth]{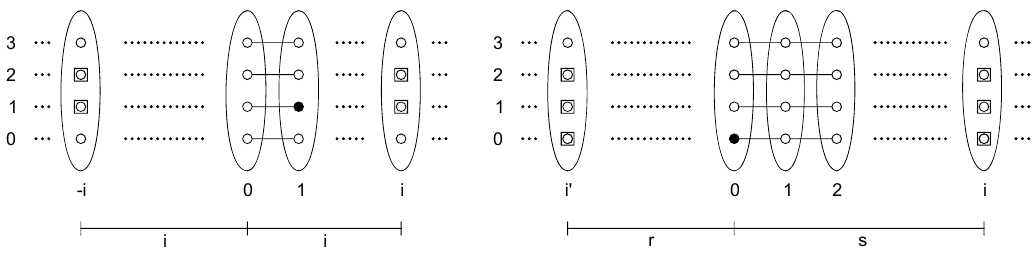}
\caption{Black vertex resolves two squared vertices  of different
columns if $i'=-i$ (left) and if $(i+i')/2\in \{0,1\}$, and
therefore $r\le s-2$ (right).}\label{fig.prop9n4}
\end{center}
\end{figure}

Observe that such a $k$ exists except for the case $n=4$ and $\{
j,j',\frac{i+i'}2 \}=\{ 0,1,2 \}$. In any other case, In any other
case, $|k-i|\not= |k-i'|$. Therefore, $d_{P_{2\infty }\Box K_n
}((i,j),(k,k))=|k-i|+1\not=|k-i'|+1= d_{P_{2\infty }\Box K_n
}((i',j'),(k,k))$. Hence, $(k,k)\in S$ resolves $x$ and $y$ (see
Figure~\ref{fig.prop9cas2i}).

It only remains to prove that $S=\{ (0,0),(1,1),(2,2)\}$ resolves
$x$ and $y$ in $P_{2\infty }\Box K_4$ when  $\{ j,j',\frac
{i+i'}2\}=\{0,1,2\}$.

Case $\frac{i+i'}2=0$. This implies $i'=-i\not=0$ and $\{j,j'
\}=\{ 1,2\}$. We may assume $i>0$. Then, $d_{P_{2\infty }\Box
K_4}((1,1),(i,j))=(i-1)+d_{K_4} (1,j)\le
(i-1)+1=i<i+1=d_{P_{2\infty }}(1,-i)=d_{P_{2\infty }}(1,i')\le
d_{P_{2\infty }}(1,i')+d_{ K_4}(1,j')=d_{P_{2\infty }\Box
K_4}((1,1),(i',j'))$. Therefore, $(1,1)$ resolves $x$ and $y$ (see
Figure~\ref{fig.prop9n4}, left).

Case $\frac{i+i'}2\in \{ 1,2 \}$. Assume $i>i'$. Then, $i\ge 2$
and $i'\le 0$, if $\frac{i+i'}2=1$, or $i'\le 1$, if
$\frac{i+i'}2=2$. For $i'\le 0$, $d_{P_{2\infty }\Box
K_4}((0,0),(i,j))=i+d_{K_4}(0,j)\ge i>  -i'+1\ge
-i'+d_{K_4}(0,j')=d_{P_{2\infty }\Box K_4}((0,0),(i',j'))$.
Therefore, $(0,0)$ resolves $x$ and $y$  (see
Figure~\ref{fig.prop9n4}, right). Finally, the case $i'=1$ is only
possible for $\frac {i+i'}2=2$, and therefore $i=3$, $\{ j, j'
\}=\{ 0,1 \}$. Then, $d_{P_{2\infty }\Box K_4}((0,0),(1,j'))=
1+d_{K_4} (0,j') \le 2<3=d_{P_{2\infty }}(0,3)\le d_{P_{2\infty
}\Box K_4}((0,0),(3,j))$. Thus, $(0,0)$ resolves $x$ and $y$.
\end{proof}

We summarize the results obtained in the preceding propositions in
Table~\ref{tab.cartProductPaths} and illustrate the given metric
basis in Figures~\ref{fig.baseCamins}, \ref{fig.baseCicles} and
\ref{fig.baseComplets}.
\begin{table}[ht]
 \begin{center}
 \begin{tabular}{|c|cccc|}
  \hline
   $G\backslash H$& \phantom{m}$P_n$, $n\ge 2$\phantom{m} & \phantom{m}$C_{n}$, $n\ge  3$ odd\phantom{m} &
    \phantom{m}$C_{n}$, $n\ge 4$ even\phantom{m} & \phantom{m}$K_n$, $n\ge 4$\phantom{m}   \\
  \hline
 $P_{\infty}$& 2&2&3& $n-1$ \\
 $P_{2\infty}$& 3&3&4
 & $n-1$\\
  \hline
 \end{tabular}
 \end{center}
\caption{Metric dimension of $G\Box H$ for some families of graphs.}
\label{tab.cartProductPaths}
\end{table}



\end{document}